\documentclass[12pt]{amsart}
\newtheorem{theorem}{Theorem}
\newtheorem{proposition}{Proposition}
\newtheorem{definition}{Definition}

 \newtheorem{corollary}{Corollary}

\usepackage[latin1]{inputenc}

\input{amssym.def}

\begin{document}

\title{Isometries for the Carath\'eodory metric} 

\author{Marco Abate and Jean-Pierre Vigu\'e}
\maketitle

\section{Introduction}

The following problem has been studied by many authors. Let $D_{1}$ and $D_{2}$ 
be two bounded domains in complex Banach spaces and let 
$ f\colon D_{1}\to D_{2} $ 
be a holomorphic map such that $ f'(a) $ is a surjective isometry for the 
Carath\'eodory infinitesimal 
metric at a point $ a $ of $ D_{1}$. The problem is to know whether 
$ f $ is an analytic isomorphism of $ D_{1} $ onto $ D_{2}$. For example, 
J.-P. Vigu\'e \cite{v1} proved this is the case when
$ D_{1} $ and $ D_{2} $ are two bounded domains in $ {\Bbb
 C}^{n} $ and $ D_{1} $ is convex. Similar results have been obtained when 
$ D_{2} $ is convex using the Kobayashi infinitesimal metric 
(I. Graham \cite{gr} and L. Belkhchicha
 \cite{be}). We have to remark that all these results are based on the theorem 
of L. Lempert (\cite{le} et \cite{le2}; one can also consult 
M. Jarnicki and P. Pflug \cite{ja}) on the equality of Kobayashi and Carath\'eodory 
metrics on a bounded convex domain in $ {\Bbb C}^{n}$. 
J.-P. Vigu\'e \cite{v2} proved the first results on this subject in the case of 
bounded domains in complex Banach spaces.

Now, we can study the same problem dropping the hypothesis that 
 $ f' (a) $ is surjective. So, we only suppose
that $ f' (a) $ is an isometry for the Carath\'eodory infinitesimal metric.
Does this imply that $ f(D_{1}) $ is a complex analytic closed 
submanifold of  $D_{2} $ and that  $ f $ is an analytic isomorphism
of $ D_{1} $ onto $ f(D_{1})$? 

Some results have been obtained by J.-P. Vigu\'e \cite{v3} and
 P. Mazet \cite{ma} assuming that $ D_{1} $ 
and $ D_{2} $ are open unit balls in complex Banach spaces, that $ a=0 $,
 and that the image of $ f' (0) $ contains enough complex extremal points 
 of the boundary of $ D_{2}$. Under these hypotheses they proved that $ f $ is
 linear equal to $ f' (0)$. This result shows that $ f(D_{1}) $ 
is an analytic submanifold of $ D_{2} $ and that $ f $ is an analytic
isomorphism 
of $ D_{1} $ onto $ f(D_{1})$. 

Of course, if we do not suppose the existence of complex extremal points in the 
image of $ f' (0)$, the map $f$ has no reason to be linear. However, one
can hope that $ f(D_{1}) $ still
is a complex analytic submanifold of $ D_{2} $. In this paper we shall be able
to prove such a result for maps of unit balls of complex Banach spaces, under some 
additional
hypotheses on the Banach spaces involved.

The authors would like to thank the referee for his/her useful remarks.

\section{The main results}

We shall prove the following theorem:

\begin{theorem}Let $(E_{j},\| \|_{j}) $ be complex Banach spaces and
 let $ B_{j}=\{x\in E_{j}\mid \|x\|_{j}<1\}$, for~$j=1,2$. Let $
f\colon B_{1}\to B_{2}
$  be a holomophic mapping with $ f(0)=0 $  and $\|f'(0)(X)\|_{2}=\|X\|_{1}$ for
all~$X\in E_{1}$.
 Then the following statements are equivalent :
\begin{enumerate}
\item there exists a direct decomposition $E_{2}=f'(0)(E_{1})\oplus F$  such that
 the  corresponding projection $\pi\colon E_{2}\to f'(0)(E_{1})$ has norm~$1$;
\item $f( B_1)$ is a closed complex direct submanifold of $B_2$, the map $f$ is a 
biholomorphism of $B_1$ onto $f(B_{1})$, and there exists a holomorphic retraction of 
$B_2$ onto $f(B_{1})$.
\end{enumerate}
\end{theorem}

To apply this theorem, we give the following definition:
 
\begin{definition} We say that a pair $(E_1,E_2)$ of complex Banach spaces 
has the property {\rm (V)} if for every linear isometry 
$L\colon E_{1}\to E_{2}$ there exists a direct decomposition $E_{2}=L(E_{1})\oplus F$
 such that the corresponding projection $\pi\colon E_{2}\to L(E_{1})$ has norm 1.
\end{definition}

From theorem 1 and definition 1, we deduce the following

\begin{theorem} Assume that the pair $(E_1,E_2)$ of complex Banach spaces has the property 
{\rm (V)},
and let $B_1$ and ${B_2}$ be their open unit balls.  Let $f\colon B_1\to {B_2}$ be a
holomorphic map such that 
\begin{enumerate}
\item $f(0)=0$, and $f'(0)$ is an isometry for the Carath\'eodory infinitesimal metric,

or 

\item $B_1$ and ${B_2}$ are homogeneous, and there exists $a\in B_1$ such that  
$f'(a)$ is an isometry for the Carath\'eodory infinitesimal metric.
\end{enumerate}
\smallskip\noindent
Then $f( B_1)$ is a closed complex direct submanifold of $B_2$, the map $f$ is a 
biholomorphism of $B_1$ onto $f(B_{1})$, and there exists a holomorphic retraction of 
$B_2$ onto $f(B_{1})$.

\end{theorem}

Now, we clearly need examples of pairs of complex Banach spaces satisfying property (V). 
The first (easy) example is given by Hilbert spaces.

\begin{proposition} Let $E_2$ a complex Hilbert space. Then, for every complex Banach 
space~$E_1$ the
pair
$(E_1,E_2)$ has property {\rm (V)}. 
\end{proposition}

More interesting is the following theorem:

\begin{theorem} Let $ I $ be a set and let $ l^{\infty}(I) $ 
be the complex Banach space of bounded sequences indexed by $I$, 
with the usual norm. Let $E_2$ be any Banach space.
 Then, the pair $\bigl( l^{\infty}(I),E_2\bigr)$ has property {\rm (V)}.
\end{theorem}

Other pairs enjoying property (V) can be constructed using suitable subspaces 
of~$\ell^\infty(I)$. For
instance, let $c_0(I)\subset\ell^\infty(I)$ be the subspace given by the elements 
$(a_i)_{i\in
I}\in\ell^\infty(I)$ such that for every $\varepsilon>0$ there exists a finite subset 
$K\subseteq I$ so
that $|a_i|<\varepsilon$ when $i\notin K$. Then:

\begin{theorem} For any sets $I$, $J$ the pair 
$\bigl(c_{0}(I),c_{0}(J)\bigr) $ has property~{\rm (V)}.
\end{theorem}

Applying Theorem~2 and~3 with $I$ finite, we get in particular a new result in the
finite-dimensional case:

\begin{corollary}
Let $f\colon\Delta^n\to D$ be a holomorphic map between a polydisk $\Delta^n\subset{\Bbb 
C}^n$
and an open convex circular bounded domain 
$D\subset{\Bbb C}^N$ (i.e., $D$ is the unit ball for a suitable norm in ${\Bbb C}^N$).
We also assume~$n\le N$, and that $D$ is homogeneous
 (for instance, $D=\Delta^N$, $B^N$ or a bounded symmetric domain). Assume that there
exists
$a\in\Delta^n$ such that $f'(a)$ is an isometry for the Carath\'eodory infinitesimal 
metrics. Then
$f(\Delta^n)$ is a closed complex submanifold of~$D$, the map $f$ is a biholomorphism onto 
its
image, and $f(\Delta^n)$ is a holomorphic retract of~$D$.
\end{corollary}

Before proving these results, we need to recall some facts.

\section{Some classical results} 

The definition and the main properties of Carath\'eodory and Kobayashi infinitesimal 
metrics $ E_{D} $ and $ F_{D} $ on a bounded domain  $ D$ are given in the book of 
T. Franzoni et E. Vesentini \cite{fr} (see also the book of M. Jarnicki and 
P. Pflug \cite{ja}). 

Let  $ B $ be the open unit ball of a complex Banach space $ E$. It is well known that 
$$
E_{B}(0,x)=F_{B}(0,x)=\|x\|\;.
$$
Furthermore, every biholomorphism $ f\colon D_{1}\to D_{2} $ 
between domains in complex Banach spaces is 
an isometry for the Carath\'eodory and Kobayashi infinitesimal metrics. 

Finally, let us recall that, the open unit balls~$B$ of the complex Banach spaces 
$c_{0}(I)$ and
 $ l^{\infty}(I) $ are
 homogeneous. Indeed, it is easy to check that, for every $ a\in
 B$, the map $\varphi_a\colon B\to B$ given by 
$$
\forall i\in I\qquad\qquad \varphi_{a}(f)_i=\frac{f_i+a_i}{1+\overline{a_i}f_i}\;,
$$
is an analytic automorphism of $B$. 

Another example of homogeneous unit ball is given by the open unit ball~$B$ of the
space $  C(S,{\Bbb C})$
of continuous complex functions on a compact space~$S$, because for every 
$a\in B$ the
map $\varphi_a\colon B\to B$ given by 
$$
\varphi_{a}(f)=\frac{f+a}{1+\overline{a}f}\;
$$
is a biholomorphism of $B$.

\section {Proof of Theorems 1 and 2} 

To begin, let us prove theorem 1.

\smallskip

{\it Proof of Theorem 1.} First, if $r:B_{2}\longrightarrow f(B_{1})$ is a holomorphic 
retraction,
$r' (0)$ is a projection of norm $\leq1$ for the Carath\'eodory infinitesimal metrics, and, 

as the 
Carath\'eodory infinitesimal metric at the origin is equal to the given norm, we get $\|r' 
(0)\|=1$. 
This proves that (2)  implies (1).

To prove that (1) implies (2), let us consider 
$$
\varphi=\pi{\circ}f\colon B_1\to f'(0)(E_1)\;.
$$
We have $\varphi(0)=0$, $\varphi(B_1)\subseteq f'(0)(E_1)\cap B_2$ (because $\pi$ has 
norm~1), and
$\varphi' (0)=\pi{\circ}f' (0)=f' (0)$.
So $\varphi' (0)$ is a linear isometry from $E_1$ onto $f' (0)(E_1)$.
Using Cartan's uniqueness theorem (see \cite{fr}), one easily proves that $\varphi$ is
a  linear isometry
from~$B_1$ onto $B'_{2} =f' (0)(E_{1})\cap B_2$.

Finally, let $\psi\colon B'_2\to F$ be defined by 
$$
\psi (y)=({\rm id}-\pi)\bigl(f(\varphi^{-1}(y))\bigr)\;.
$$ 
Then the set $f(B_{1})$ is the graph of $\psi$, the map $(\pi,\psi\circ\pi)\colon B_2\to 
f(B_1)$ is a
holomorphic retraction of~$B_2$ onto~$f(B_1)$, and $\varphi^{-1}\circ\pi|_{f(B_1)}\colon 
f(B_1)\to B_1$
is a holomorphic inverse of~$f$, and the theorem is proved.\hfill\textsquare
\medskip

Now, we can prove Theorem 2.
\medskip

{\it Proof of Theorem 2.}  First, let us remark that, in case (2), 
by pre-composing  $ f $ with an analytic automorphism of $ B_{1} $ and post-composing it 
with
an analytic automorphism of $ B_{2}$, we can assume that 
$ f(0)=0 $ and that $ 0 $ is precisely the point $a$ such that
$f' (0)$ is an isometry for the Carath\'eodory infinitesimal metrics. Thus without loss of 
generality in
both cases we can assume that $f' (0)$ is an isometry for the norms of~$E_1$ and $E_2$. 

Since $(E_1,E_2)$ satisfies the property (V), there exists a direct decomposition
$E_{2}=f'(0)(E_{1})\oplus F$ such that the corresponding projection $\pi\colon E_{2}\to
f'(0)(E_{1})$ has norm  1 and we can apply Theorem 1.\hfill\textsquare

\section {Examples of pair of Banach spaces with property (V)}

Now we have to give examples of pair of complex Banach spaces satisfying property (V). 
Proposition
1  (the case of Hilbert spaces) is easy and left as an exercise. Let us now give the
\bigskip

{\it Proof of Theorem 3.} We suppose that $E_1=l^{\infty}(I)$  and we consider an 
isometry $L\colon \ell^\infty(I) \to E_2$. Let 
$ G\colon L(E_{1}) \to l^{\infty}(I)$  be the inverse of $L$.
So, $G$ is a linear map of norm 1; for every $i\in I$, let $G_i$
be the $i$-component of $G$. Then $G_i$ is a linear form from $L(E_{1})$ to~${\Bbb C}$
of norm 1. By the Hahn-Banach Theorem, we can extend $G_i$
to a linear form $H_i\colon E_2\to{\Bbb C}$ still of norm 1. Setting
$H=(H_{i})_{i\in I}$ we obtain a linear map 
$H\colon E_{2} \to l^{\infty}(I)$ of norm 1 extending~$G$.
Then it is clear that $L{\circ}H$ is a projection of~$E_2$ onto $L(l^{\infty}(I))$ of 
norm~1, and
taking $F=\hbox{Ker}(L\circ H)$ the theorem is proved.\hfill\textsquare

\bigskip

{\it Proof of Theorem 4.} 
Let $ L\colon c_{0}(I)\to c_{0}(J)$ be an isometry, and let 
$ (e^k)_{k\in I}$ be the
canonical basis of $ c_{0}(I)$. Since $L$ is an isometry, for every $ k\in
I$ there exists $ j(k)\in J$ such that $ |L(e^{k})_{j(k)}|=1$.
 Now, if we consider an element $ v=(v_{i})_{i\in I} $ of $ c_{0}(I)$ 
such that $ v_{k}=0$, then $ L(v)_{j(k)}=0$. In fact, suppose that
 $ L(v)_{j(k)}\neq0$. For $ \lambda\in{\Bbb C} $ small enough, we have
  $ \|e^{k}+\lambda v\|=1$. But 
$$
L(e^{k}+\lambda v)_{j(k)}=L(e^{k})_{j(k)}+\lambda L(v)_{j(k)}
=e^{i\theta}+\lambda
 L(v)_{j(k)}\;.
$$
Therefore if $ L(v)_{j(k)}\neq0$, there exists $ \lambda\in{\Bbb C} $ small enough 
such that the modulus of $ L(e^k+\lambda v)_{j(k)}$ is greater than 1, and thus 
$\|L(e^k+\lambda
v)\|>1$, contradiction. It follows that the map $k\mapsto j(k)$ is injective.

Let $M=\{j(k)\mid k\in I\}\subseteq J$, and let 
$ \pi\colon c_{0}(J)\to c_{0}(M) $ be the canonical projection. The previous argument shows 

that
$\pi\circ L(e^k)=\lambda_k e^{j(k)}$ with $|\lambda_k|=1$ for all~$k\in I$; it is then easy 

to check
that $\varphi=\pi{\circ}L\colon c_{0}(I)\to c_{0}(M)$ is a linear surjective isometry, and
that $L\circ\varphi^{-1}\circ\pi\colon c_0(J)\to L\bigl(c_0(I)\bigr)$ is a linear 
projection of norm~1
of~$c_0(J)$ onto~$L\bigl(c_0(I)\bigr)$, as required.\hfill\textsquare

\bigskip

It might be interesting to remark that the same proof yields that a pair $(E_1,E_2)$ of 
complex Banach
spaces satisfies property (V) if each $E_j$ has a Schauder basis $(e^k_j)$ such that
$$
\left\|\sum_k \lambda_k e^k_j\right\|_{E_j}=\sup_{k}|\lambda_k|\;.
$$  
 
\section {Final remarks}

Not all pairs of complex Banach spaces have property (V); so we do not know whether 
Theorem 2 holds in general.

For example, the Banach spaces  $ c_{0}(\Bbb N)$ is not complemented in 
$l^{\infty}(\Bbb N)$, and so the pair $\bigl(c_{0}(\Bbb N)$, $l^{\infty}(\Bbb N)\bigr)$ 
does not have property (V). 

It is also possible to build finite dimensional examples. 
Take $E_2=({\Bbb C}^3,\|\cdot\|_\infty)$, so that the unit ball of~$E_2$ is the open
polydisk~$\Delta^{3}$. If $L\colon{\Bbb C}^2\to{\Bbb C}^3$ is given by $L(x,y)=(x,y,x+y)$,
then $B=L^{-1}(\Delta^3)$ is the open unit ball in $ {\Bbb C}^{2} $ for a norm~$\|\cdot\|$;
set $E_1=({\Bbb C}^2,\|\cdot\|)$. We claim that the pair $(E_1,E_2)$ does not satisfy 
property (V). By
construction, $L\colon E_1\to E_2$ is a linear isometry. The set $
(1,0,1)+(\{0\}\times\Delta\times\{0\})$ is contained in the boundary of
$\Delta^3$. Since
$(1,0,1)\in L(\partial B)$, it is easy to check that if there exists a projection~$\pi$
 of norm 1 from $ {\Bbb C}^{3} $ onto $ L({\Bbb C}^{2})$, then~$\pi$
 must vanish on $ \{0\}\times{\Bbb C}\times\{0\}$. Considering the point
 $ (0,1,1)$, we analogously see that $\pi$ must vanish on $\Delta\times\{0\}\times\{0\}$; 
and thus
such a projection cannot exist.

\bigskip\bigskip

Marco Abate, Dipartimento di Matematica, Universit\`a di Pisa,Largo Pontecorvo 5, 56127 
Pisa, Italy.
\smallskip

e-mail :  abate@dm.unipi.it

\medskip

J.-P. V., LMA, Universit\'{e} de Poitiers, CNRS, 
Math\'{e}matiques, SP2MI, BP 30179, 86962 FUTUROSCOPE.
\smallskip

e-mail :  vigue@math.univ-poitiers.fr

\end{document}